\documentclass[12pt]{amsart}
\usepackage{amsmath,amscd,amssymb,amsfonts}
\def\ver{Jan.~10, 2008, v.6}

\def\ssbull{\raise.2ex\hbox{${\scriptscriptstyle\bullet}$}}
\def\scirc{\,\raise.2ex\hbox{${\scriptstyle\circ}$}\,}
\def\msum{\hbox{$\sum$}}
\def\mprod{\hbox{$\prod$}}
\def\mopl{\hbox{$\bigoplus$}}
\def\d{{\rm d}}
\def\fm{{\mathfrak m}}

\def\bC{{\mathbf C}}
\def\bN{{\mathbf N}}
\def\bP{{\mathbf P}}
\def\bQ{{\mathbf Q}}
\def\bR{{\mathbf R}}
\def\bZ{{\mathbf Z}}
\def\cD{{\mathcal D}}
\def\cH{{\mathcal H}}
\def\cI{{\mathcal I}}

\def\cO{{\mathcal O}}
\def\cP{{\mathcal P}}

\def\ophi{\bar{\phi}}
\def\talpha{\widetilde{\alpha}}
\def\rd{\partial}
\def\DR{{\rm DR}}
\def\Gr{{\rm Gr}}
\def\IH{{\rm IH}}
\def\Sing{{\rm Sing}\,}
\def\Spec{{\rm Spec}\,}
\def\Coker{{\rm Coker}}
\def\Ker{{\rm Ker}}
\def\Im{\hbox{\rm Im}}
\def\prim{{\rm prim}}
\def\({{\rm (}}
\def\){{\rm )}}

\begin{document}
\title[Griffiths theorem on rational integrals, II]
{A generalization of Griffiths theorem\\
on rational integrals, II}
\author{Alexandru Dimca }
\address{Laboratoire J.A. Dieudonn\'e, UMR du CNRS 6621,
Universit\'e de Nice-Sophia Antipolis, Parc Valrose,
06108 Nice Cedex 02, FRANCE.}
\email{dimca@math.unice.fr}
\author{Morihiko Saito }
\address{RIMS Kyoto University, Kyoto 606--8502 JAPAN}
\email{ msaito@kurims.kyoto-u.ac.jp}
\author{Lorenz Wotzlaw }
\address{Fachbereich Mathematik und Informatik II,
Mathematisches Institut, Freie Universit\"at Berlin,
Arnimallee 3, D-14195 Berlin, Germany}
\email{wotzlaw@math.fu-berlin.de}
\date{\ver}
\begin{abstract}
We show that the Hodge and pole order filtrations are globally
different for sufficiently general singular projective hypersurfaces
in case the degree is 3 or 4 assuming the dimension of the projective
space is at least 5 or 3 respectively.
We then study an algebraic formula for the global Hodge filtration in
the ordinary double point case conjectured by the third named author.
This is more explicit and easier to calculate than the previous one
in this case.
We prove a variant of it under the assumption that the image of the
singular points by the $e$-fold Veronese embedding consists of
linearly independent points, where $e$ is determined only by the
dimension and the degree.
In particular, the original conjecture is true in case the above
condition is satisfied for $e=1$.

\end{abstract}
\maketitle

\bigskip\centerline{\bf Introduction}

\bigskip\noindent
Let $X=\bP^n$, and $Y\subset X$ be a hypersurface defined by a reduced
polynomial $f$ of degree $d$. Set $U=X\setminus Y$.
Let $F,P$ denote respectively the global Hodge and pole order
filtrations on the cohomology $H^n(U,\bC)$, see [5], [6].
Locally it is easy to calculate the difference between these two
filtrations at least in the case of isolated weighted homogeneous
singularities, see (1.3.2) below.
However, this is quite nontrivial globally (i.e.\ on the cohomology).
It is important to know when the two filtrations coincide globally,
since the Hodge filtration and especially the Kodaira-Spencer map
can be calculated rather easily if they coincide, see [9], Thm.~4.5.
It is known that they are different if $Y$ has bad singularities
(see [7] and also [9], 2.5).
In case the singularities consist of ordinary double points,
however, it was unclear whether they still differ globally.
They coincide for $n=2$ in this case (loc.~cit.)
but the calculation for the case $n>2$ is quite complicated in general.
In this paper we show

\medskip\noindent
{\bf Theorem~1.} {\it
Assume $d=3$ with $n\ge 5$ or $d=4$ with $n\ge 3$.
Set $m=[n/2]$, and assume $1+(n+1)/d\le p \le n-m$.
Then, for a sufficiently general singular hypersurface $Y$, we have
$F^p\ne P^p$ on $H^n(U,\bC)$.
}

\medskip
Here a {\it sufficiently general} singular hypersurface means
that it corresponds to a point of a certain (sufficiently small)
non-empty Zariski-open subset of $D\setminus\Sing D$, where $D$ is
the parameter space of singular hypersurfaces of degree $d$ in
$\bP^n$, see (3.6).
In particular, $\Sing Y$ consists of one ordinary double point.
It is unclear whether the two filtrations differ whenever
$\Sing Y$ consists of one ordinary double point.
According to Theorem~1, the formula for the Kodaira-Spencer map
in [9], Thm.~4.5 is effective only for $p>n-m$ in the ordinary
double point case.
By Theorem~2 below, however, we can show a similar formula in the
ordinary point case which is valid also for $p\le n-m$, see (4.5).
In the $n$ odd case, we can also use the self-duality for the
calculation of the Kodaira-Spencer map, see Remark~(3.9)(ii).

Theorem~1 implies $F^{p}\ne P^{p}$ on $H^{n+1}_Y(X,\bC)$
by the long exact sequence associated with local cohomology.
If $n$ is odd and $Y$ has only ordinary double points as
singularities, then $Y$ is a $\bQ$-homology manifold so that
$H^{n-1}(Y,\bQ)$ coincides with the intersection
cohomology $\IH^{n-1}(Y,\bQ)$ (see [2], [10]),
and also with the local cohomology $H^{n+1}_Y(X,\bQ)(1)$.
In particular, they have a pure Hodge structure in this case.
If $n=3$, we cannot calculate directly $F^1$ on $\IH^2(Y,\bQ)$, but
this can be obtained from $F^2$ if we can calculate the
intersection pairing.
For example, if $(n,d)=(3,4)$, then $Y$ is a singular $K3$ surface,
i.e. its blow-up along the singular points is a smooth $K3$ surface,
and there is a lot of work on the lattice and the intersection
pairing.

Let $R=\bC[x_0,\dots,x_n]$ with $x_0,\dots,x_n$ the coordinates of
$\bC^{n+1}$.
Let $J\subset R$ be the Jacobian ideal of $f$ (i.e. generated by
$f_j:=\rd f/\rd x_j$), and
$I$ be the ideal generated by homogeneous functions
vanishing at the singular points of $Y$.
Let $R_k$ denote the degree $k$ part of $R$, and similarly for
$I_k$, etc. Set $q=n-p$, $m=[n/2]$, and $I^j=R$ for $j\le 0$.
Assume that $\Sing\,Y$ consists of ordinary double points.
Then the third named author ([23], 6.5) found the following

\medskip\noindent
{\bf Conjecture 1.}
$\,\,\Gr_F^{p}H^n(U,\bC)=(I^{q-m+1}/I^{q-m}J)_{(q+1)d-n-1}$.

\medskip\noindent
This is a generalization of Griffiths theorem on rational integrals
[12], but it is quite different from the one in [9].
Indeed, the formula in [9], Th.~1 is for the case of general
singularities, and it is not necessarily easy to calculate concrete
examples because of the problems of torsion and inductive limit
which produce a problem of infinite dimensional vector spaces
making explicit calculations quite difficult.
Conjecture~1 is much more explicit and algebraic
(or ring theoretic).
It is much easier to calculate concrete examples by using Conjecture~1
in the ordinary double point case.
For the moment the relation between the two generalizations
of the Griffiths theorem is unclear, since the results of [9]
imply only that $\Gr_F^{p}H^n(U,\bC)$ is a quotient of
$(I^{q-m+1}/fI^{q-m})_{(q+1)d-n-1}$.

The original argument in [23] was essentially correct for $p\ge n-m$
(using [9], [19], [20]).
Actually Conjecture~1 holds for such $p$ in the case of general
singularities by modifying $m$ and $I$ appropriately, see (2.2).
In the case $p<n-m$, however, there are some difficulties:
among others, the coincidence of the Hodge and pole order filtrations,
which is not true as is shown in Theorem~1, was used
(in fact, this problem was rather extensively studied there using
the theory of logarithmic forms for strongly quasi-homogeneous
singularities, see e.g. a remark after Thm.~3.14 in [23]).
For other difficulties, see (2.3.1), (2.3.4) below.

Let $\cI\subset\cO_X$ be the reduced ideal of $\Sing\,Y\subset X$. Set
$I^{(i)}_k=\Gamma(X,\cI^i(k))$ and $I^{(i)}=\mopl_k\,I^{(i)}_k$.
The difference between $I^i$ and $I^{(i)}$ is one of the main
problems, see remarks after (2.3.1).
We have by definition the exact sequences
$$
0 \longrightarrow I^{(i)}_k\longrightarrow R_k
\buildrel{\beta_{k}^{(i)}}\over\longrightarrow \mopl_{y\in\Sing Y}\,
\cO_{X,y}/\fm_{X,y}^{i},
\leqno(0.1)
$$
choosing a trivialization of $\cO_{X,y}(k)$, where $\fm_{X,y}=\cI_y$
is the maximal ideal of $\cO_{X,y}$.
In this paper we prove a variant of Conjecture~1 as follows:

\medskip\noindent
{\bf Theorem~2.} {\it
Assume the singular points are ordinary double points.
For $q=n-p>m=[n/2]$, we have canonical isomorphisms
$$
\aligned
\Gr_F^{p}H^n(U,\bC)&=(I^{(q-m+1)}/I^{(q-m)}J)_{(q+1)d-n-1}\\
&=(I^{(q-m+2)}/(I^{(q-m+2)}\cap I^{(q-m)}J))_{(q+1)d-n-1},
\endaligned
\leqno(0.2)
$$
if the following condition is satisfied in the notation of
{\rm (0.1):}

\begin{itemize}
\item[($A$)]
$\beta_{k}^{(i)}$ is surjective for
$(k,i)=(qd-n,q-m+1)$ and $(qd-n-1,q-m)$.
\end{itemize}

\noindent
Moreover, condition~$(A)$ is satisfied if we have the
following\,{\rm :}

\begin{itemize}
\item[($B$)]
For $e=m(d-1)-p$, the image of the singular points by the
$e$-fold Veronese embedding consists of linearly independent
points.
\end{itemize}
}

\medskip
Note that $(I^{(q-m)}J)_{(q+1)d-n-1}=\msum_{j=0}^n\,
f_jI^{(q-m)}_{(q+1)d-n-d}$.
Condition~$(B)$ means that, for each singular point $y$, there is
a hypersurface of degree $e$ containing the singular points
other than $y$, but not $y$,
see (2.3.5) below.
In order to satisfy $(B)$, there should hold at least the inequality
$|\Sing Y|\le\binom{e+n}{n}$.
For $n$ even, this is always satisfied by Varchenko~[22].
For $n$ odd, however, this is not necessarily satisfied, e.g. if
$n=3$, $d=4$, $q=2$, and $Y$ is a Kummer surface with
$16$ ordinary double points where condition~$(A)$ is not satisfied
either but Conjecture~1 seems to hold.
There seem to be some examples such that condition~$(A)$ is satisfied
but $(B)$ is not, see (4.7) below.
The proof of Theorem~2 uses the theory of Brieskorn modules [3]
in the ordinary double point case by restricting to a neighborhood
of each singular point, see (4.1--3).

In a special case, we can deduce from Theorem~2 and (2.5) below
the following:

\medskip\noindent
{\bf Corollary~1.} {\it
Conjecture {\rm 1} is true if the singular points consist
of ordinary double points and are linearly independent points
in $\bP^n$ \(in particular, $|\Sing\,Y|\le n+1)$.
}

\medskip
In general Conjecture 1 is still open.

\medskip
In Sect.~1, we review some basic facts from the theory of Hodge
and pole order filtrations for a hypersurface of a smooth variety.
In Sect.~2, we study the case of hypersurfaces of projective spaces.
In Sect.~3, we prove Theorem~1 by constructing examples explicitly.
In Sect.~4, we prove Theorem~2 and Corollary~1 after reviewing some
basic facts about Brieskorn modules in the ordinary double point case.

\newpage
\centerline{\bf 1. Hodge and pole order filtrations}

\bigskip\noindent
{\bf 1.1.}
Let $X$ be a proper smooth complex algebraic variety of dimension
$n\ge 2$, and $Y$ be a reduced divisor on $X$. Set $U=X\setminus Y$.
Let $\cO_X(*Y)$ be the localization of the structure
sheaf $\cO_X$ along $Y$.
We have the Hodge filtration $F$ on $\cO_X(*Y)$.
This is uniquely determined by using the relation with the
$V$-filtration of Kashiwara~[14] and Malgrange~[16], see [17].
Moreover, $F$ induces the Hodge filtration $F^p$ of
$H^j(U,\bC)$ by taking the $j$-th cohomology group of the subcomplex
$F^p\DR(\cO_X(*Y))$ defined by
$$
F_{-p}\cO_X(*Y)\to\cdots\to F_{n-p}\cO_X(*Y)\otimes\Omega_X^n.
\leqno(1.1.1)
$$
Indeed, this is reduced to the normal crossing case by using a
resolution of singularities together with the stability of mixed
Hodge modules by the direct image under a proper morphism.
In this case the Hodge filtration $F$ on $\cO_X(*Y)$ is given by
using the sum of the pole orders along the irreducible components,
and the assertion follows from [4] as is well known.

Let $P$ be the pole order filtration on $\cO_X(*Y)$ (see [6]), i.e.
$$
P_i\cO_X(*Y)=\cO_X((i+1)Y)\,\,\,\text{for $i\ge 0$, and
$0$ otherwise.}
$$
Note that the pole order filtration in [4], II, (3.12.2) is by the sum
of the orders of poles along the irreducible components in the normal
crossing case, and actually coincides with our Hodge filtration $F$
on the de Rham complex.

If $Y$ is smooth, then $F_i=P_i$ on $\cO_X(*Y)$
(see also [11], [12]). So we get in the general case
$$
F_i\subset P_i\,\,\,\text{on}\,\,\,\cO_X(*Y).
$$
Let $h$ be a local defining equation of $Y$ at $y\in Y$,
$b_{h,y}(s)$ be the $b$-function of $h$, and
$\talpha_{Y,y}$ be the smallest root of $b_{h,y}(-s)/(1-s)$.
Then we have by [18]
$$
F_i=P_i\,\,\,\text{on $\cO_{X,y}(*Y)$ if $i\le \talpha_{Y,y}-1$}.
\leqno(1.1.2)
$$
If $y$ is an ordinary double point, then
$b_{h,y}(s)=(s+1)(s+n/2)$ and hence $\talpha_{Y,y}=n/2$
as is well known.
Note that (1.1.2) was first obtained by Deligne at least if
$h$ is a homogenous polynomial of degree $r$ with an isolated
singularity (where $\talpha_{Y,y}=n/r$),
see e.g. Remark 4.6 in [18].

As a corollary of (1.1.2), we get
$$
\hbox{$F^p=P^p$ on $H^j(U,\bC)$ if $p\ge j-\talpha_Y+1$},
\leqno(1.1.3)
$$
where
$\talpha_Y=\min\{\talpha_{Y,y}\,|\,y\in\Sing Y\}$.
Indeed, $P$ on $H^j(U,\bC)$ is defined by the image of the
$j$-th cohomology group of the complex $P^p\DR(\cO_X(*Y))$ as in
(1.1.1) with $F$ replaced by $P$, and this coincides with the image
of the cohomology group of the subcomplex
$\sigma_{\le j}P^p\DR(\cO_X(*Y))$,
where $\sigma_{\le j}$ is the filtration ``b\^ete" in [5], i.e.
$$
\sigma_{\le j}P^p\DR(\cO_X(*Y))=
[P_{-p}\cO_X(*Y)\to\cdots\to P_{j-p}\cO_X(*Y)\otimes\Omega_X^j].
\leqno(1.1.4)
$$
Indeed, the $k$-th cohomology group of the quotient complex of
$P^p\DR(\cO_X(*Y))$ by (1.1.4) vanishes for $k\le j$.

In the case $j=n=3$ or $4$ and $\Sing Y$ consists of ordinary double
points as in Theorem~1, we have $\talpha_Y=n/2$, $n-m=2$, and
the equality in (1.1.3) holds for $p\ne n-m$.

\medskip\noindent
{\bf 1.2.~Local cohomology.}
Since $H^j(X,\DR(\cO_X(*Y)/\cO_X))=H^{j+1}_Y(X,\bC)$,
we get the Hodge and pole order filtrations on $H^{j+1}_Y(X,\bC)$ in a
similar way.
Moreover, we have the compatibility of the long exact sequence
$$
\to H^j(X,\bC)\to H^j(U,\bC)\to H^{j+1}_Y(X,\bC)\to H^{j+1}(X,\bC)\to,
\leqno(1.2.1)
$$
with the pole order filtration (i.e. it is exact after taking $P^p$)
if $X=\bP^n$.

Indeed, we have a short exact sequence
$$
0\to P_i\cO_X\to P_i\cO_X(*Y)\to P_i(\cO_X(*Y)/\cO_X)\to 0,
\leqno(1.2.2)
$$
where the filtration $P$ on $\cO_X$ and $\cO_X(*Y)/\cO_X$ are
respectively the induced and quotient filtrations.
This induces the long exact sequence
$$
H^j(P^p\DR(\cO_X))\buildrel{\alpha_j}\over\to
H^j(P^p\DR(\cO_X(*Y)))\buildrel{\beta_j}\over\to
H^j(P^p\DR(\cO_X(*Y)/\cO_X)),
$$
where the cohomology group is taken over $X$, and the filtration
$P$ on $\DR(\cO_X))$ and $\DR(\cO_X(*Y)/\cO_X))$ are defined as in
(1.1.1) with $F$ replaced by $P$.
Since $X=\bP^n$, the restriction morphism $H^j(X,\bC)\to H^j(U,\bC)$
vanishes for $j\ne 0$, and the long exact sequence splits into
a family of short exact sequences.
This implies $\alpha_j=0$ for $j\ne 0$ using $F=P$ on $\DR(\cO_X)$
because $\alpha_j$ with $P$ replaced by $F$ vanishes by the strictness
of the Hodge filtration $F$ on $\bR\Gamma(X,\DR(\cO_X(*Y)))$.
So the assertion follows from the snake lemma using the
strictness of $F=P$ on $\bR\Gamma(X,\DR(\cO_X))$.

\medskip\noindent
{\bf 1.3.~Semi-weighted-homogeneous case.}
Assume $Y$ has only isolated singularities which are locally
semi-weighted-homogeneous, i.e. $Y$ is analytically locally defined
by a holomorphic function $h = \sum_{\alpha\ge 1}h_{\alpha}$,
where the $h_{\alpha}$ for $\alpha\in\bQ$ are weighted homogeneous
polynomials of degree $\alpha$ with respect to some local coordinates
$x_{1},\dots,x_{n}$ around $y\in\Sing Y $ and some positive weights
$w_{1},\dots,w_{n}$, and moreover, $h_{1}^{-1}(0)$ (and hence $Y$)
has an isolated singularity at $y$.
In this case, it is well known that
$$
\talpha_{Y,y} = \msum_{i}\, w_{i},
\leqno(1.3.1)
$$
by Kashiwara's unpublished work (this also follows from [15] together
with [3]).

Let $\cO_{X,y}^{\ge \beta}$ be the ideal of $\cO_{X,y}$ generated by
$\prod_{i}x_{i}^{\nu_{i}}$ with $\sum_{i}w_{i}\nu_{i}\ge\beta-
\talpha_{Y,y}$.
Let $\cD_{X}$ be the sheaf of linear differential operators with the
filtration $F$ by the order of differential operators.
Put $k_{0} = [n-\talpha_{Y,y}]-1$. Then we have by [19]
$$
\aligned
F_{p}(\cO_{X,y}(*Y))
&=\msum_{k\ge 0}\,F_{p-k}\cD_{X,y}(\cO_{X,y}^{\ge k+1}h^{-k-1})\\
&=\msum_{k=0}^{k_{0}}\,F_{p-k}\cD_{X,y}(\cO_{X,y}^{\ge k+1}h^{-k-1}).
\endaligned
\leqno(1.3.2)
$$
If $w_i=1/b$ for any $i$ with $b\in\bN$, then (1.3.2) implies for
$p=m:=[\talpha_{Y,y}]$
$$
F_{m}(\cO_{X,y}(*Y))=\cO_{X,y}^{\ge m+1}h^{-m-1}.
\leqno(1.3.3)
$$
This does not hold in general, e.g. in case the weights
are $\frac{1}{3},\frac{1}{3},\frac{1}{2},$ with $n=3$.

\medskip\noindent
{\bf 1.4.~Ordinary double point case.}
Assume that $\Sing Y$ consists of ordinary double points.
Then $b_{h,y}=(s+1)(s+n/2)$ and hence $\talpha_{Y,y}=n/2$
as is well known (see also (1.3.1)). Set $m=[n/2]$.
Then $k_0=m-1$ and $\cO_{X,y}^{\ge k+1}=\cO_{X,y}$ for $k\le k_0$.
So (1.3.2) becomes
$$
F_{p}(\cO_{X,y}(*Y))=
F_{p-m+1}\cD_{X,y}(\cO_{X,y}h^{-m})\,\,\,\,\text{if}\,\, p\ge m-1,
\leqno(1.4.1)
$$
where $F_{p}(\cO_{X,y}(*Y))=P_p(\cO_{X,y}(*Y))$ if $p<m-1$.

This implies the following lemma, which is compatible with (1.1.2),
and was conjectured by the third named author (see [20]):

\medskip\noindent
{\bf 1.5.~Lemma.} {\it With the above notation and assumption,
we have
$$
F_p(\cO_X(*Y))=\cI^{p-m+1}\cO_X((p+1)Y)\,\,\,\hbox{for}\,\,\,p\ge 0,
\leqno(1.5.1)
$$
where $\cI$ is the reduced ideal of $\Sing Y\subset X$, and
$\cI^{p-m+1}=\cO_X$ for $p\le m-1$.
}

\medskip\noindent
{\it Proof.}
We reproduce here an argument in [20].
By (1.4.1) it is enough to show the following by increasing induction
on $p\ge 0$:
$$F_p\cD_{X,y}h^{-m}=\cI_y^ph^{-m-p}.
\leqno(1.5.2)
$$
Here $\cI$ is the maximal ideal at $y$, and we may assume
$h=\sum_{i=1}^nx_i^2$ using GAGA if necessary.
We have to show by increasing induction on $p\ge 0$
$$
u=x^{\nu}h^{-m-p}\in F_p\cD_{X,y}h^{-m}\,\,\, \text{if}\,\, |\nu|=p,
\leqno(1.5.3)
$$
where $x^{\nu}=\prod_ix_i^{\nu_i}$ for $\nu= (\nu_1,\dots,\nu_n)\in\bN^n$.
Here we may assume $\nu_i\ne 1$ for any $i$ and $p>1$, because the
assertion is easy otherwise.
Then we have $x^{\nu}=x_i^2x^{\mu}$ for some $i$, and
$$\partial_i(x_ix^{\mu}h^{-(m+p-1)})=((\mu_i+1)h-(m+p-1)x_ih_i)x^{\mu}
h^{-m-p}.$$
Adding this over $i$, we get (1.5.3), because $|\mu|+n-2(m+p-1)\ne 0$.
So (1.5.2) and hence (1.5.1) follow.

\bigskip\bigskip
\centerline{\bf 2. Projective hypersurface case}

\bigskip\noindent
{\bf 2.1.~Hodge filtration.}
With the notation of (1.1), assume $X=\bP^n$ with $n\ge 2$.
Then we have by [9], Prop.~2.2
$$
H^k(X,F_p\cO_X(*Y))=0\quad\hbox{for}\,\,\,k>0.
\leqno(2.1.1)
$$
As a corollary, $F^pH^j(U,\bC)$ is given by the $j$-th cohomology of
the complex
$$
\Gamma(X,F_{-p}\cO_X(*Y))\to\cdots\to\Gamma(X,F_{n-p}\cO_X(*Y)
\otimes\Omega_X^n).
$$

Let $R=\bC[x_0,\dots,x_n]$, where
$x_0,\dots,x_n$ are the coordinates of $\bC^{n+1}$.
Let $J$ be the ideal of $R$ generated by
$f_i:=\rd f/\rd x_i\,(0\le i\le n)$.
Let $R_k$ denote the degree $k$ part of $R$ so that
$R=\mopl_k R_k$, and similarly for $J_k$, etc.
Let
$$
\xi=\hbox{$\frac{1}{d}$}\msum_i\, x_i\rd/\rd x_i,
$$
so that $\xi f=f$.
Let $\iota_\xi$ denote the interior product by $\xi$.
Let $\Omega^j$ be the vector space of global algebraic (i.e.
polynomial) $j$-forms on $\bC^{n+1}$, and let $\Omega^j[f^{-1}]_k$
be the degree $k$ part of $\Omega^j[f^{-1}]$,
where the degrees of $x_i$ and $\d x_i$ are 1.
Then
$$\iota_{\xi}(\Omega^{j+1}[f^{-1}]_0)=\Gamma(U,\Omega_U^j).
\leqno(2.1.2).$$
This is compatible with the differential $\d$ up to a sign,
because
$$
\iota_{\xi}\scirc\d+\d\scirc\iota_{\xi}=L_{\xi},
$$
with $L_{\xi}$ the Lie derivation and $L_{\xi}\eta=(k/d)\eta$ for
$\eta\in(\Omega^j[f^{-1}])_k$.
We have for $g\in R$
$$
\d(gf^{-k}\omega_i)=(-1)^i(f\partial_i g-k g f_i)f^{-k-1}\omega,
\leqno(2.1.3)
$$
where
$\omega=\d x_0\wedge\cdots\wedge\d x_n$,
$\omega_i=\d x_0\wedge\cdots\widehat{\d x_i}\cdots\wedge\d x_n$.

Let
$$
m=[\talpha_Y].
$$
For $q\in\bN$, let $\cI_{(q)}$ be the ideal of $\cO_X$ such that
$$
F_q(\cO_X(*Y))=\cI_{(q)}\cO_X((q+1)Y).
\leqno(2.1.4)
$$
Then $\cI_{(q)}=\cO_X$ for $q<m$ by (1.1.2).
Let
$$
I_k=\Gamma(X,\cI_{(m)}(k))\subset R_k,\quad I=\mopl_{k\in\bN}I_k\subset
R.
$$
Taking local coordinates $y_0,\dots,y_n$ of $\bC^{n+1}\setminus\{0\}$
such that $\rd/\rd y_0=\xi$, we get
$$
\iota_{\xi}(I\Omega^{n+1})=\Im\,\iota_{\xi}\cap I\Omega^n,
\leqno(2.1.5)
$$
using the injectivity of
$$
\iota_{\xi}:\Omega^{n+1}[f^{-1}]\to\Omega^n[f^{-1}].
$$
Here we have also the following argument:
For $g\in R_k$, we have $g\in I_k$ if and only if $x_ig\in I_{k+1}$
for any $i\in[0,n]$. (This follows from the definition of $I$.)

Note that $m=\talpha_Y=+\infty$ if $Y$ is smooth, and
$\talpha_{Y,y}=\sum_{i=0}^n w_i$ in case $Y$ is
analytically locally defined by a semi-weighted-homogeneous
function $h$ with weights $w_0,\dots,w_n$ at $y\in\Sing\,Y$,
see (1.3.1).

\medskip
From (2.1.1--5) we can deduce a generalization of a theorem of
Griffiths~[12] as follows
(here no condition on the singularities of $Y$ is assumed).

\medskip\noindent
{\bf 2.2.~Theorem.} {\it With the above notation \(e.g.
$m=[\talpha_Y])$, we have
}
$$
\Gr_F^{n-q}H^n(U,\bC)=\begin{cases}
(R/J)_{(q+1)d-n-1}&\text{if}\,\,\,q<m,\\
(I/J)_{(q+1)d-n-1}&\text{if}\,\,\,q=m.
\end{cases}
\leqno(2.2.1)
$$

\medskip\noindent
{\it Proof.}
Since $f\in J$, the assertion immediately follows from (2.1.1--5).

\medskip\noindent
{\bf 2.3.~Ordinary double point case.}
Assume $\Sing Y$ consists of ordinary double points so that
$m=[\talpha_Y]=[n/2]$ as in (1.4).
Then $\cI_{(m)}$ in (2.1.4) coincides with the (reduced) ideal $\cI$
of $\Sing Y\subset X$ by (1.5.1).
Without the assumption on the singularities, this does not hold,
see (1.3.3).
Using (2.1.1) and (1.5.1), the third named author obtained (2.2.1)
in this case (i.e. Conjecture~1 for $p\ge n-m$), see [23], 6.5.

Let
$$
I^{(i)}_k=\Gamma(\bP^{n},\cI^i(k))\subset R_k,\quad
I^{(i)}=\mopl_k\,I^{(i)}_k\subset R.
$$
Then $(I^i)_k\subset I^{(i)}_k\subset R_k$, but it is not clear whether
$$
(I^i)_k=I^{(i)}_k.
\leqno(2.3.1)
$$
Note that (2.3.1) holds for $k\gg 0$, because the restriction to
$\Spec R\setminus \{0\}$ of the sheaf corresponding to $I^i$ coincides
with that for $I^{(i)}$.
However, (2.3.1) for an arbitrary $k$ does not hold in general if
$q\ge 2$.
For example, let $f=xyz(x+y+z)$ with $n=2$.
In this case there is no hypersurface of degree $\le 2$ passing through
all the six singular points of $Y$, i.e. $I_i=0$ for $i\le 2$.
So $g\in I^{(2)}_4\ne (I^2)_4=0$. See also (2.4) below.

Choosing a section of $\cO_X(1)$ not vanishing at $y\in\Sing Y$,
we can trivialize $\cO_{X,y}(k)$ so that we get exact sequences
$$
0\longrightarrow I^{(i+1)}_{k} \longrightarrow I^{(i)}_{k}
\buildrel{\gamma_{k}^{(i)}}\over \longrightarrow
\mopl_{y\in\Sing Y}\,\fm_{X,y}^{i}/\fm_{X,y}^{i+1},
\leqno(2.3.2)
$$
where $\fm_{X,y}=\cI_y$ is the maximal ideal of $\cO_{X,y}$. Let
$$
I^{(i),(y)}_{k}=\Ker(\gamma_{k}^{(i)}:I^{(i)}_{k}\to
\mopl_{y'\in\Sing Y\setminus\{y\}}\,\fm_{X,y'}^{i}/\fm_{X,y'}^{i+1}).
$$
If $\gamma_k^{(i)}$ is surjective, then we have the surjectivity of
$$
\gamma_{k}^{(i),(y)}:I^{(i),(y)}_{k}\to\fm_{X,y}^{i}/\fm_{X,y}^{i+1},
\leqno(2.3.3)
$$
where $\gamma_{k}^{(i),(y)}$ is the restriction of $\gamma_{k}^{(i)}$.

By (1.5.1) and (2.1.2), we get an injection
$$
\iota_{\xi}\big((I^{(j-p-m+1)}\Omega^{j+1})_{(j-p+1)d}\,
f^{-(j-p+1)}\big)\hookrightarrow
\Gamma(U,F_{j-p}\cO_X(*Y)\otimes\Omega_X^j).
$$
Here
$(I^{(i)}\Omega^{j})_{k}=I^{(i)}_{k-j}\otimes_{\bC}(\Omega^{j})_{j}$,
because $\Omega^j=R\otimes_{\bC}(\Omega^j)_j$.

One of the main problems is whether the above injection is
surjective, i.e.
$$
\iota_{\xi}(I^{(i')}\Omega^{j+1})_{k'}
=\Im\,\iota_{\xi}\cap(I^{(i')}\Omega^j)_{k'},
\leqno(2.3.4)
$$
where $i'=j-p-m+1$, $k'=(j-p+1)d$.
Note that (2.3.4) for $j=n$ holds by the same argument as in the proof
of (2.1.5).
However, (2.3.4) for $j<n$ does not hold, for example, if $i'=k'-j$
(without assuming that $i',k'$ are as above).

In (2.6--8) below, we will show that (2.3.4) is closely related to
the surjectivity of (2.3.3) and also to the following:

\begin{itemize}
\item[(2.3.5)]
For each $y\in\Sing Y$, there is $g_{(y)}\in\Gamma(X,\cO_X(e))$
such that $y\notin g_{(y)}^{-1}(0)$ and $\Sing Y\setminus\{y\}\subset
g_{(y)}^{-1}(0)$, where $e$ is a given positive integer.
\end{itemize}
This condition is satisfied for any $e'>e$ if it is satisfied
for $e$.
(Indeed, it is enough to replace $g_{(y)}$ with $h_{(y)} g_{(y)}$
where $h_{(y)}$ is any section of $\cO_X(e'-e)$ such that $y\notin
h_{(y)}^{-1}(0)$.)
Condition~(2.3.5) means that the images of the singular points by the
$e$-fold Veronese embedding $i_{(e)}$ in (3.6) correspond to linearly
independent vectors in the affine space.

\medskip\noindent
{\bf 2.4.~Linearly independent case.}
Assume the singular points correspond to linearly independent vectors
in $\bC^{n+1}$.
Replacing the coordinates if necessary, we may assume
$\Sing Y=\{P_0,\dots,P_s\}$ where $s\in[0,n]$ and the $P_i$
are defined by the $i$-th unit vector of $\bC^{n+1}$.
In this case $I^{(i)}\subset R$ is a monomial ideal, and for a monomial
$x^{\nu}:=\prod_jx_j^{\nu_j}$ we have
$$
x^{\nu}\in I^{(i)}\Leftrightarrow \text{$x^{\nu}|_{x_j=1}\in
\fm_j^i$ for each $j\in[0,s]$,}
\leqno(2.4.1)
$$
where $\fm_j$ is the maximal ideal generated by
$x_l\,(l\ne j)$.
Let $\Gamma^{(i)}\subset\bN^{n+1}$ such that
$$
I^{(i)}=\msum_{\nu\in\Gamma^{(i)}}\, \bC x^{\nu}.
$$
Set $|\nu|_{(j)}=\msum_{k\ne j}\nu_k$. Then we have
$$
\Gamma^{(i)}=\big\{\nu\in\bN^{n+1}\,\big|\,|\nu|_{(j)}\ge i\,\,
(j\in[0,s])\big\}.
\leqno(2.4.2)
$$
If $|\nu|=k$, then the condition $|\nu|_{(j)}\ge i$ is
equivalent to $\nu_j\le k-i$.
If $i=1$, then $I$ is generated by $x_j$ for $j>s$ and
$x_jx_l$ for $j,l\in[0,s]$ with $j\ne l$.

In the case $s=0$, we have
$$
I^{(i)}=I^i\quad\text{for any}\,\, i\ge 1\,\,\,\text{if}\,\,\,
|\Sing Y|=1.
\leqno(2.4.3)
$$
Assume $s=n$ for simplicity.
Then $I$ is generated by $x_ix_j$ for $i\ne j$,
and $I^{(2)}$ is generated by $x_i^2x_j^2$ for $i\ne j$
and $x_ix_jx_l$ for $i,j,l$ mutually different.
So we get $I^{(2)}_k=(I^2)_k$ for $k\ge 4$, but
$I^{(2)}_3\ne(I^2)_3=0$.

\medskip
More generally, we have the following:

\medskip\noindent
{\bf 2.5.~Lemma.} {\it 
Assume the singular points of $Y$ correspond to linearly independent
vectors in $\bC^{n+1}$. Then
}
$$
(I^i)_k=I^{(i)}_k\quad\text{if}\,\,\,k\ge 2i.
\leqno(2.5.1)
$$

\medskip\noindent
{\it Proof.}
We may assume $i\ge 2$ and $s\ne 0$ by (2.4.3).
With the notation of (2.4), any $x^{\nu}\in I^{(i)}_k$ is divisible
either by $x_j$ with $j>s$ or by $x_jx_l$ with $j,l\in[0,s]\,(j\ne l)$.
(Indeed, otherwise $x^{\nu}=x_j^k$ for some $j\in[0,s]$,
but $x_j^k\notin I^{(i)}$.)
So we can proceed by increasing induction on $i$, applying the
inductive hypothesis to the case where $i,k$ are replaced by
$i-1$ and $k-2$ respectively.

\medskip\noindent
{\bf 2.6.~Lemma.} {\it
Assume $\Sing Y$ consists of ordinary double points, and
{\rm (2.3.5)} is satisfied for $e=k-i(d-1)$.
Then $\gamma_{k}^{(i)}$ in {\rm (2.3.2)} is surjective so that
we get a short exact sequence
$$
0\longrightarrow I^{(i+1)}_{k} \longrightarrow I^{(i)}_{k}
\buildrel{\gamma_{k}^{(i)}}\over \longrightarrow
\mopl_{y\in\Sing Y}\,\fm_{X,y}^{i}/\fm_{X,y}^{i+1} \longrightarrow 0,
\leqno(2.6.1)
$$
where $\fm_{X,y}=\cI_y$ is the maximal ideal of $\cO_{X,y}$.
}

\medskip\noindent
{\it Proof.}
For each $y\in\Sing Y$, the $f_j\in I_{d-1}$ for $j\in[0,n]$
generate $\cI_y=\fm_{X,y}$, and hence the $g_{(y)}\prod_j f_j^{\nu_j}$
for $|\nu|=i$ generate $\fm_{X,y}^{i}/\fm_{X,y}^{i+1}$.
So the assertion follows.

\medskip\noindent
{\bf 2.7.~Remarks.} (i)
The morphism $\beta_{k}^{(j)}$ in (0.1) is surjective if and only if
$\gamma_{k}^{(i)}$ in (2.3.2) is surjective for any $i\in[0,j-1]$.
So Lemma~(2.6) shows that condition~$(B)$ in Theorem~2 implies $(A)$,
since $qd-n-(q-m)(d-1)=m(d-1)-p$ and $d\ge 2$.

\medskip
(ii) Let $g=\sum_{|\nu|=k}a_{\nu}x^{\nu}\in R_k$.
Then $g\in I^{(i)}_k$ if and only if
$$
\text{$(\rd^{\mu}g)(y)=0$ for any $y\in\Sing Y$ and $\mu\in\bN^{n+1}$
with $|\mu|=i-1$,}
\leqno(2.7.1)
$$
where $\rd^{\mu}g=\prod_{i=0}^n\rd_i^{\mu_i}g$.
Let $M=\binom{i-1+n}{n}|\Sing Y|$ and $N=\binom{k+n}{n}$.
The $a_{\nu}$ are viewed as coordinates of $\bC^N$
parametrizing the homogeneous polynomials of degree $k$,
and (2.7.1) gives $M$ linear relations among the $a_{\nu}$ defining
the subspace $I_k^{(i)}\subset R_k$.
So $\beta_k^i$ is surjective if and only if these $M$ relations are
linearly independent, i.e. the corresponding matrix of size
$(M,N)$ has rank $M$.

\medskip\noindent
{\bf 2.8.~Proposition.} {\it
Assume $\Sing Y$ consists of ordinary double points.
Then {\rm (2.3.4)} with $j=n-1$ holds if $\gamma_{k}^{(i)}$ in
{\rm (2.3.2)} is surjective for $k=k'-n-1$ and any $i\in[0,i'-1]$.
}

\medskip\noindent
{\it Proof.}
By increasing filtration on $i'>0$, it is enough to show
$$
\iota_{\xi}(\eta)\in\iota_{\xi}(I^{(i')}\Omega_n)_{k'}
\,\,\,\text{if}\,\,\,
\eta\in (I^{(i'-1)}\Omega_n)_{k'}\,\,\,\text{with}\,\,\,
\iota_{\xi}(\eta)\in (I^{(i')}\Omega^{n-1})_{k'}.
\leqno(2.8.1)
$$
For each $y\in\Sing\,Y$, take coordinates $x^{(y)}_0,\dots,x^{(y)}_n$
such that $y=(1,0,\dots,0)$.
With the notation of Lemma~(2.6), set $k=k'-n-1$.
Then in the notation of (2.3.3), the hypothesis of the proposition
implies the surjectivity of
$$
\gamma_{k}^{(i'-1),(y)}:I^{(i'-1),(y)}_{k}\to
\fm_{X,y}^{i'-1}/\fm_{X,y}^{i'}.
$$
So we may replace $\eta$ with
$\sum_{y}x_0^{(y)}\eta^{(y)}$ where
$$
\eta^{(y)}\in I^{(i'-1),(y)}_{k}\otimes_{\bC}(\Omega^n)_n
\quad\text{with}\quad
\gamma_{k+1}^{(i'-1)}(\eta)=\msum_{y}\,
\gamma_{k+1}^{(i'-1),(y)}(x_0^{(y)}\eta^{(y)}).
$$
Then, for the proof of (2.8.1), we may assume
$$
\eta\in x_0^{(y)}I^{(i'-1),(y)}_{k}\otimes_{\bC}
(\Omega^n)_n\,\,\,\,\text{for some}\,\,\, y\in\Sing Y,
$$
because for any $g\in I^{(i'-1),(y)}_{k}$, $g|_{X\setminus\{y\}}$
is a section of $\cI^{(i')}(k)|_{X\setminus\{y\}}$.

Let
$\omega^{(y)}=\d x^{(y)}_0\wedge\cdots\wedge\d x^{(y)}_n$,
$\omega^{(y)}_j=\d x^{(y)}_0\wedge\cdots\widehat{\d x^{(y)}_j}\cdots
\wedge\d x^{(y)}_n$.
Then
$$
\eta=\msum_{j=0}^n\,x_0^{(y)}h_j^{(y)}\omega^{(y)}_j\quad
\text{with}\,\,\, h_j^{(y)}\in I^{(i'-1),(y)}_{k}.
$$
Calculating mod $I^{(i')}\Omega^{n-1}$ the coefficient of
$$
\d x^{(y)}_1\wedge\cdots\widehat{\d x^{(y)}_j}\cdots\wedge
\d x^{(y)}_n \,\,\,\,\text{in}\,\,\,\, \iota_{\xi}
(\msum_{j=0}^n\,x^{(y)}_0 h^{(y)}_j\omega^{(y)}_j),
$$
which belongs to $I^{(i')}\Omega^{n-1}$ by the hypothesis of (2.8.1),
we see that $h_j^{(y)}\in I^{(i')}_{k}$ for $j\ne 0$.
Then we may assume $h_j^{(y)}=0$ for $j\ne 0$ so that
$$
\eta=x^{(y)}_0 h^{(y)}_0\omega^{(y)}_0.
$$
By the definition of $I^{(i'-1),(y)}_{k}$, we have
$$
x^{(y)}_jh^{(y)}_0\omega^{(y)}_j\in I^{(i')}\Omega^n\quad
\text{for}\,\,\,j\ne 0,
$$
and
$$
\msum_{j=0}^n\,(-1)^j\iota_{\xi}(x^{(y)}_jh^{(y)}_0\omega^{(y)}_j)
=\iota_{\xi}(\iota_{\xi}(h^{(y)}_0\omega^{(y)}))=0.
$$
So the assertion follows.

\newpage
\centerline{\bf 3. Proof of Theorem~1}

\bigskip\noindent
{\bf 3.1.~Problem.}
Assume $X=\bP^n$ and $\Sing Y$ consists of ordinary double points.
One of the main problems in generalizing a theorem of Griffiths~[12]
is whether the following equality holds:
$$
\aligned
&\qquad\qquad F^pH^n(U,\bC)=P^pH^n(U,\bC),\\ \hbox{i.e.}
\quad&\Im\bigl(\iota_{\xi}\bigl((\Omega^{n+1})_{(q+1)d}
\,f^{-(q+1)}\bigr)\to H^n(U,\bC)\bigr)\subset F^pH^n(U,\bC),
\endaligned
\leqno(3.1.1)
$$
where $q=n-p$.
This was rather extensively studied in [23]
(see e.g. a remark after Th. 3.14 there).
We show that (3.1.1) does not hold in general, see (3.7--8).
This implies that the isomorphism in Conjecture~1 for $p<n-m$
(i.e. $q>m$) cannot be deduced by the method indicated there.

\medskip\noindent
{\bf 3.2.~Proposition.} {\it
Let $X,Y$ be as above. Assume $q=n-p>m$ and
$F^{p+1}=P^{p+1}$ on $H^{n}(U,\bC)$.
Then $\Gr_F^pH^n(U,\bC)$ is a subquotient of $(I/J)_{(q+1)d-n-1}$.
}

\medskip\noindent
{\it Proof.}
By (2.1), $H^n(U,\bC)$ is the cokernel of
$$
\d:\Gamma(X,\Omega^{n-1}_{X}(*Y))\to\Gamma(X,\Omega^n_{X}(*Y)),
$$
and $P^pH^n(U,\bC)$ is the image of
$\Gamma(X,(\cO_X((q+1)Y))\otimes\Omega^n_{X})$, and similarly for
$F$.

Let $\cI$ be the reduced ideal of $\Sing Y\subset X$, and $I_k=
\Gamma(X,\cI(k))\subset R_k$.
By assumption together with Lemma~(1.5) we have
$$
F_q(\cO_Z(*Y))\subset \cI\cO_X((q+1)Y),\,\,\,\text{and}\,\,\,
F^{p+1}=P^{p+1}\,\,\,\text{on}\,\,\,H^{n}(U,\bC).
\leqno(3.2.1)
$$
So we get a commutative diagram
$$
\CD
\Gamma(X,F_{q-1}\cO_X\otimes\Omega_X^{n-1})\oplus
\Gamma(X,F_{q-1}\cO_X\otimes\Omega_X^n) @>>>
\Gamma(X,F_{q}\cO_X\otimes\Omega_X^n),\\
@VVV @VVV\\
\Gamma(X,\Omega_X^{n-1}(qY))\oplus\Gamma(X,\Omega_X^n(qY))
 @>{\phi}>> \Gamma(X,\cI\Omega_X^n((q+1)Y)),\\
\endCD
$$
By (2.1.2--3) together with the inclusion $Rf\subset J$, we have
$$
\Coker\,\phi=(I/J)_{(q+1)d-n-1}.
\leqno(3.2.2)
$$
So the assertion is reduced to
$$
\text{$\Gr_F^pH^n(U,\bC)$ is a subquotient of $\Coker\,\phi$.}
\leqno(3.2.3)
$$
Taking the image of the diagram by the canonical morphism to
$H^n(U,\bC)$ and adding the cokernels, we get
$$
\CD
F^{p+1}H^n(U,\bC) @>>> F^{p}H^n(U,\bC) @>>>
\Gr_F^{p}H^n(U,\bC) @>>> 0,\\
@| @VV{\cap}V @VV{\cap}V\\
P^{p+1}H^n(U,\bC) @>{\ophi}>> P^{p}H^n(U,\bC) @>>>
\Coker\,\ophi @>>> 0,\\
\endCD
$$
where the image of $\d\Gamma(X,\Omega_X^{n-1}(qY))$ in $H^n(U,\bC)$
vanishes (considering the case $q=\infty$).
Moreover, $\Coker\,\ophi$ is a quotient of $\Coker\,\phi$
by the snake lemma.
So the assertion follows.

\medskip\noindent
{\bf 3.3.~Hodge numbers of smooth hypersurfaces.}
Define integers $C(n+1,d,i)$ by
$$
(t+\cdots+t^{d-1})^{n+1}=\msum_{i=n+1}^{(n+1)(d-1)}\,C(n+1,d,i)t^i,
\leqno(3.3.1)
$$
so that
$$
C(n+1,d,i)=C(n+1,d,(n+1)d-i),
$$
where $C(n+1,d,i)=0$ unless $i\in[n+1,(n+1)(d-1)]$.
This is the Poincare polynomial of the graded vector space
$$
\Omega^{n+1}/dg\wedge\Omega^n,
$$
if $g$ is a homogeneous polynomial
of degree $d$ with an isolated singularity at the origin
(e.g. if $g=\sum_ix_i^d$).
For the hypersurface $Z'\subset X$ defined by $g$,
we have by Griffiths~[12]
$$
C(n+1,d,pd)=\dim \Gr_F^{n-p}H^{n-1}_{\prim}(Z',\bC)\quad\text{for}
\,\,\, p\in [1,n],
\leqno(3.3.2)
$$
where $H^{n-1}_{\prim}(Z',\bC)$ denotes the primitive part.

\medskip\noindent
{\bf 3.4.~Isolated singularity case.}
Assume $Y$ has only isolated singularities, and
$$
n-p>q_0:=\max\{q\,|\,\Gr_F^qH^{n-1}(F_y,\bC)\ne 0\,\,\,\text{for some}
\,\, y\in\Sing Y\},
\leqno(3.4.1)
$$
where $F$ is the Hodge filtration on the vanishing cohomology
$H^{n-1}(F_y,\bC)$ at $y\in\Sing Y$, see [21].
Here $F_y$ denotes the Milnor fiber around $y$.
In the case $\Sing Y$ consist of ordinary double points, we have
$q_0=m:=[n/2]$, see (3.5.1) below.
In general, we have $q_0\ge (n-1)/2$ by the Hodge symmetry.

Under the above assumptions we have
$$
\dim\Gr_F^pH^n(U,\bC)=\dim\Gr_F^pH_Y^{n+1}(X,\bC)=C(n+1,d,pd).
\leqno(3.4.2)
$$
Indeed, there is a perfect pairing of mixed Hodge structures
$$
H_Y^{n+1}(X,\bQ)\times H^{n-1}(Y,\bQ)\to\bQ(-n),
\leqno(3.4.3)
$$
and condition~(3.4.1) (together with $q_0\ge (n-1)/2$) and (3.3.2)
imply
$$
\dim\Gr_F^{n-p}H^{n-1}(Y,\bC)=C(n+1,d,pd).
\leqno(3.4.4)
$$
The last assertion is reduced to the case of a smooth hypersurface
by taking a 1-parameter deformation
$Z_t=:\{f+tg=0\}\,(t\in \Delta)$ of $Y=Z_0$ whose general fibers $Z_t$
and total space $Z$ are smooth where we assume that the hypersurface
$\{g=0\}$ does not meet $\Sing Y$.
Here we use also the exact sequence of mixed Hodge structures
$$
0\to H^{n-1}(Y)\to H^{n-1}(Z_{\infty})\buildrel{\rho}\over\to
\mopl H^{n-1}(F_y)\to H^{n}(Y)\to H^{n}(Z_{\infty})\to 0,
\leqno(3.4.5)
$$
(see also [8], 1.9), where $H^{n-1}(Z_{\infty})$ denotes the limit
mixed Hodge structure.
Note that $\Gr_F^{n-p}H^{n-1}(Z_{\infty},\bC)=\Gr_F^{n-p}
H^{n-1}_{\prim}(Z_{\infty},\bC)$ because $n-p>(n-1)/2$.

\medskip\noindent
{\bf 3.5.~Remark.}
Assume the singularities of $Y$ are ordinary double points.
Since the weight filtration on the unipotent (resp. non-unipotent)
monodromy part of $H^{n-1}(F_y,\bQ)$ has the symmetry with center $n$
(resp. $n-1$) by definition~[21], and the monodromy on the vanishing
cycles is $(-1)^n$, we have
$$
H^{n-1}(F_y,\bQ)=\bQ(-m),
\leqno(3.5.1)
$$
where $m=[n/2]$.
In particular, $\rho$ in (3.4.5) is surjective for $n$ odd
(considering the monodromy), and we get by the above argument
$$
|\Sing\,Y|\le C(n+1,d,(m+1)d) \quad\text{if}\,\,\, n=2m+1.
\leqno(3.5.2)
$$
This is related to [22].
Note that $\rho$ can be non-surjective if $n$ is even and
$\Sing Y$ consists of sufficiently many ordinary double points.
Indeed, the Betti number $b_n(Y)$ may depend on the position of the
singularities, see for instance Thm. (4.5) in [7], p. 208.
The position of singularities enters there via the dimension of
$I_{md-2m-1}$ (where $n=2m$).
Its proof uses an exact sequence
$$
P^{m+1}H^n(\bP^n\setminus Y)\to \mopl_{y\in\Sing Y}
H^n(B_y\setminus Y)\to H_0^n(Y)(-1)\to 0,
\leqno(3.5.3)
$$
where $B_y\subset\bP^n$ is a sufficiently small ball with center
$y$. Here $H_0^n(Y)$ denotes the primitive cohomology defined by
$\Coker(H^n(\bP^n)\to H^n(Y))$.
Note that $$H^n(B_y\setminus Y)=\Coker(N:H^{n-1}(F_y)\to
H^{n-1}(F_y)(-1))=\bQ(-m-1).$$

Using (3.4.5), (3.4.3), (1.2.1) and (2.2.1), we have also
$$
\aligned
\dim\Ker\,\rho&=\dim\Gr_F^{m}H^{n-1}(Y,\bC)=
\dim\Gr_F^{m}H^{n+1}_Y(X,\bC)\\
&=\dim\Gr_F^{m}H^n(U,\bC)=\dim\,(I/J)_{(m+1)d-2m-1}.
\endaligned
$$
In the case $n$ is even and $d=2$, we have
$H^{n-1}_{\prim}(Z_{\infty},\bC)=0$ and $\rho$ vanishes.
In the case $n$ is even, $d\ge 3$, and $\Sing Y$ consists of one
ordinary double point, we have $H^{n-1}_{\prim}(Z_{\infty},\bC)\ne 0$
and $\rho$ is surjective (because $\rho$ is nonzero by the theory
of vanishing cycles), and hence $b_n(Y)=1$
(for more general singularities, see Thm.~(4.17) in [7], p.214),
thus at the level of topology nothing surprising may occur.

\medskip\noindent
{\bf 3.6.~Discriminant.}
Let $i_{(d)}:X=\bP^n\to \cP=\bP^N$ be the $d$-fold Veronese embedding
defined by the line bundle $\cO_X(d)$
(i.e. by using the monomials $x^{\nu}$ of degree $d$),
where $N=\binom{n+d}{n}-1$.
Let $\cP^{\vee}$ be the dual projective space of $\cP$ parametrizing
the hyperplanes of $\cP$.
Let $\cH\subset\cP\times\cP^{\vee}$ be the universal hyperplane
whose intersection with $\cP\times\{z\}$ is the hyperplane
corresponding to $z\in\cP^{\vee}$.
Let $D$ be the discriminant of the projection
$$
pr:(i_{(d)}(X)\times\cP^{\vee})\cap\cH\to\cP^{\vee}.
$$
This is called the dual variety of $X\subset\cP$.
It is well known that $D$ is irreducible (since $D$ is the image of
a $\bP^{N-n-1}$-bundle over $X$ corresponding to the hyperplanes
which are tangent to $X$).
By the theory of Lefschetz pencils, it is also known that $\Sing Y$
consists of one ordinary double point if and only if it corresponds to a
smooth point of $D$.

\medskip\noindent
{\bf 3.7.~Proof of Theorem~1.}
By (3.6) it is enough to show
$$
F^{p+1}\ne P^{p+1}\,\,\,\text{on}\,\,\,H^n(X\setminus Y,\bC)
$$
for one hypersurface $Y$ whose singularities consist of one ordinary
double point, assuming
$$
(n+1)/d\le p < n-m,\,\,\,\text{i.e.}\,\,\, m<q\le n-(n+1)/d.
$$
Indeed, $F^{-\infty}/F^{p+1}$ defines a vector bundle on the parameter
space of hypersurfaces $Y$ whose singularities consist of one
ordinary double point, and in the notation of (3.2), $g f^{-q}\omega$
for $g\in R_{qd-n-1}$ defines a section of this bundle when $f$
varies.
Since $P^{p+1}$ is generated by these sections where $q=n-p$,
the subset defined by the condition $P^{p+1}/F^{p+1}\ne 0$ is a
Zariski-open subset.

Let
$$
f=\msum_{i=1}^n\, x_i^d/d-x_0^{d-2}\msum_{i=1}^n\, x_i^2/2,
$$
so that
$$
\aligned
f_0=-\hbox{$\frac{1}{2}$}\msum_{i=1}^3\, x_i^2,\,\,\,
f_i=x_i^{2}-x_0x_i\,(1\le i\le 3)\,\,\,\,\text{if}\,\,d=3,
\\
f_0=-\msum_{i=1}^4\, x_0x_i^2,\,\,\,
f_i=x_i^{3}-x_0^{2}x_i\,(1\le i\le 4)\,\,\,\,\text{if}\,\,d=4.
\endaligned
$$
Here $I$ is generated by $x_1,\dots,x_n$ so that $R/I=\bC[x_0]$.
By assumption, (3.2.1) and (3.4.1) are satisfied (in particular,
$q>n/2>p$), and moreover
$$
C(n+1,d,pd)\ne 0,
$$
see (3.3.1) for $C(n+1,d,k)$.
The assumptions imply also
$$
\text{$p\ge 2$, $n\ge 5$ if $d=3$, and $p\ge 1$, $n\ge 3 $ if $d=4$.}
$$

Since $q>n/2$, we get
$$
r:=(q+1)d-n-1 > d.
$$
We will show
$$
\dim\,(I/J)_r<C(n+1,d,pd),
\leqno(3.7.1)
$$
contradicting Proposition~(3.2) and (3.4.2).

Take $x^{\nu}=\mprod_{i=0}^n x_i^{\nu_i}\in I_{r}$ with
$\nu=(\nu_0,\dots,\nu_n)\in\bN^{n+1}$ where
$$
|\nu|:=\msum_{i=0}^n\, \nu_i=r,\quad \nu_0<r.
$$
Using $f_i$ for $i>0$, we can replace $x^{\nu}$ with $x^{\mu}$ mod $J_r$
(i.e. $x^{\nu}-x^{\mu}\in J_r$) so that
$$
\mu_i\le d-2\,\,(i>0),\,\,\, \nu_i-\mu_i\in(d-2)\bZ.
$$
So we may assume $\nu_i\le d-2$ for $i>0$.
Let $|\nu|'=\msum_{i=1}^n\, \nu_i$, and
$$
s=\min\{s\in\bZ\,\big|\,|\nu|'-s\in(d-2)\bZ,\,\,\, s\ge r-(d-2)\}.
$$
We first show that if $|\nu|'<r-(d-2)$ (i.e. if $\nu_0>d-2$), then
$$
x^{\nu}=(-1)^{(|\nu|'-s)/(d-2)}\msum_{|\mu|'=s,\,\mu_i\le d-2\,(i>0)}\,
e_{\nu,\mu}x^{\mu}\,\,\,\text{mod}\,\,\, J_r,
\leqno(3.7.2)
$$
where the summation is taken over $\mu$ such that $|\mu|'=s$ and
$\mu_i\le d-2$ for $i>0$, and the $e_{\nu,\mu}$ are nonnegative numbers
with $e_{\nu,\mu}\ne 0$
for some $\mu$ (for each $\nu$).
By decreasing induction on $|\nu|'$, it is enough to show (3.7.2)
with the summation taken over $b$ such that $|\mu|'=|\nu|'+(d-2)$
instead of $|\mu|'=s$.
But this modified assertion follows from
$$
\hbox{$x^{\nu}x_0^{-2}\msum_{i=1}^n\, x_i^2\in J_r\,\,$ if $\nu_0>d-2$},
\leqno(3.7.3)
$$
because we have for $i>0$ (using $f_i$)
$$
\hbox{$x^{\nu}=x^{\nu}x_0^{-2}x_i^2$ \,mod $J_r\,\,$
if $\nu_i>0$, $\nu_0\ge 2$.}
\leqno(3.7.4)
$$
(For the last argument we need the assumption $d=3$ or 4.)

Let $V_r$ be the vector space with basis $x^{\mu}$ such that
$|\mu|=r$ and $\mu_i\le d-2$ for $i\ge 0$.
Let $V_{r,k}$ be the vector subspace of $V_r$ generated by $x^{\mu}$
such that $\mu_0=k$ (i.e. $|\mu|'=r-k$).
The above argument implies that $(I/J)_r$ is spanned by
$V_r=\sum_{k=0}^{d-2}V_{r,k}$, and moreover,
$x_0^{r-2}\msum_{i=1}^n\, x_i^2\in J_r$
gives a nontrivial relation in $V_{r,r-s}$.
Thus we get (3.7.1), i.e.
$$
\dim \,(I/J)_r<\dim V_r=C(n+1,d,(q+1)d)=C(n+1,d,pd).
$$
So the assertion follows.

\medskip\noindent
{\bf 3.8.~Other examples.} (i)
It is not easy to extend the above argument to the case $d\ge 5$.
Let $n=4$, $d=5$, and
$$
f=x_0^3(x_1x_4+x_2x_3)-\msum_{i=1}^4\,x_i^5/5,
$$
so that $f_0=3x_0^2(x_1x_4+x_2x_3)$,
$f_i=x_0^3x_{5-i}-x_i^4\,\,(1\le i\le 4)$.
In this case we have
$F^2H^4(U,{\bf C})\ne P^2H^4(U,{\bf C})$ for this hypersurface
and hence for a sufficiently general singular hypersurface.

\medskip
(ii)
In the above examples, $\Sing Y$ consists of one point.
Let $n=3$, $d=4$, and
$$
f=\msum_{0\le i<j\le 3}\,x_i^2x_j^2/2,\quad
f_i=x_i\msum_{k\ne i}\, x_k^2.
$$
Then $\Sing Y$ consists of 4 points corresponding to the unit
vectors of $\bC^4$.
In this case we have
$F^2H^3(U,{\bf C})\ne P^2H^3(U,{\bf C})$ for this hypersurface.

\medskip\noindent
{\bf 3.9.~Remarks.} (i) In [9], Thm.~4.5, two of the authors gave
a formula for the Kodaira-Spencer map
$$
\Gr_F\nabla_{\xi}:\Gr_F^{p+1}H^n(U_s,\bC)\to\Gr_F^{p}H^n(U_s,\bC),
\leqno(3.9.1)
$$
where $\{Y_s\}$ is an equisingular family of hypersurfaces,
see loc.~cit.
In case the $Y_s$ have only ordinary double points,
Theorem~1 implies that the formula is useful only for $p>n-m$.
In this case, however, (3.9.1) is given by the multiplication by
$-(n-p)(\xi f)_s$ for any $p$ under the isomorphisms of Theorem~2
and Theorem~(2.2), see Theorem~(4.5) below.

\medskip
(ii) In case $n$ is odd, $Y$ is a $\bQ$-homology manifold so that
$$
H^n(U_s,\bC)=H^{n+1}_{Y_s}(X,\bC)_{\prim}=
H^{n-1}_{\prim}({Y_s},\bC)(-1),
$$
and the Kodaira-Spencer map for $p\le n-m$ can be calculated by
using the duality because the horizontality of the canonical
pairing on  $H^{n-1}_{\prim}({Y_s},\bC)$ implies
that the Kodaira-Spencer map is self-dual up to a sign.

\newpage
\centerline{\bf 4. Proof of Theorem~2}

\bigskip\noindent
{\bf 4.1.~Brieskorn modules for ordinary double points.}
We first review some basic facts about algebraic Brieskorn modules.
Let $z_1,\dots,z_n$ be the coordinates of $Z=\bC^n$, and
$h=\sum_{i=1}^n z_i^2$.
We denote by $(\Omega_Z^{\ssbull},\d)$ the complex of algebraic
differential forms on $Z$.
Let $(A_h^{\ssbull},\d)$ be the subcomplex defined by
$$
A_h^i=\Ker(\d h\wedge:\Omega_Z^i\to\Omega_Z^{i+1}).
$$
Since $(\Omega_Z^{\ssbull},\d h\wedge)$ is the Koszul complex
associated to the regular sequence $h_i=2z_i$ for $i\in[0,n]$, we have
$$
H^i(\Omega_Z^{\ssbull},\d h\wedge)=0\,\,\,\,\text{for}\,\,\,i\ne n.
\leqno(4.1.1)
$$
This implies that the cohomology group $H^iA_h^{\ssbull}$ is a left
$\bC[t]\langle\rd_t\rangle$-module for $i\ne n$, and the action of
$\rd_t^{-1}$ is well defined on the algebraic Brieskorn module
$$
H^nA_h^{\ssbull}=\Omega_{Z}^n/\d h\wedge\d\Omega_{Z}^{n-2}.
$$
Here $\rd_t[\eta]=[\phi]$ for $\eta,\phi\in A_h^i$
if there is $\sigma\in A_h^{i-1}$ such that
$$
[\eta]=[\d h\wedge\sigma],\quad [\phi]=[\d\sigma],
\leqno(4.1.2)
$$
where $[\eta]$ denotes the class of $\eta$ in
$H^iA_h^{\ssbull}$, see [3].
The action of $t$ is defined by the multiplication by $h$.
We have the finiteness of $H^iA_h^{\ssbull}$ over $\bC[t]$
by using the canonical compactification of the morphism $h$.
(The argument is essentially the same as in the analytic case in
loc.~cit.)
Then $H^iA_h^{\ssbull}$ is $t$-torsion free for $i<n$, and
we get by the theory of Milnor fibration
$$
H^iA_h^{\ssbull}=0\,\,\,\text{for}\,\,i\ne 1, n.
\leqno(4.1.3)
$$
We have the graded structure such that $\deg z_i=\deg \d z_i=1$.
This is compatible with $\d$ and $\d h\wedge$ (up to a shift of degree),
and defines a graded structure on $H^nA_h^{\ssbull}$.
Let $H^nA_{h,k}^{\ssbull}$ denote the degree $k$ part of
$H^nA_h^{\ssbull}$ so that
$$
H^nA_h^{\ssbull}=\mopl_{k\ge n}\,H^nA_{h,k}^{\ssbull}.
$$
Using the relation $\sum_i z_ih_i=2h$, we get a well-known formula
$$
2t\rd_t[\phi]=(k-2)[\phi]\,\,\,\,\text{for}\,\,\,[\phi]
\in H^nA_{h,k}^{\ssbull}.
\leqno(4.1.4)
$$
This implies the $t$-torsion-freeness of $H^nA_h^{\ssbull}$ (because
we may assume $k\ge n$).

For $i=1$, $H^1A_{h}^{\ssbull}$ is a free $\bC[t]$-module of rank one
generated by $[\d h]$. Since $A_h^0=0$, this implies
$$
H^1A_{h}^{\ssbull}=\bC[h]\d h=\Ker(\d:A_h^1\to A_h^2).
\leqno(4.1.5)
$$

Define $D'_q:\Omega_Z^i\to\Omega_Z^{i+1}$ for $q\in\bZ$ by
$$
D'_q\eta=h\d\eta-q\d h\wedge\eta.
$$
This is compatible with the graded structure up to the shift by
$\deg h=2$, and we have $D'_q\scirc D'_{q-1}=0$.
We will denote by $\Omega_{Z,j}^i$ the degree $j$ part of
$\Omega_Z^i$.

\medskip
The following will be used in the proof of Theorem~2 in (4.3).

\medskip
\newpage
\noindent
{\bf 4.2.~Lemma.} {\it Assume $j\ne 2q\ne 0$. Then

\smallskip\noindent
{\rm (i)} $D'_q:\Omega_{Z,j}^{n-1}\to\Omega_{Z,j+2}^n$ is surjective
if $j+2>n$.

\smallskip\noindent
{\rm (ii)} $\Im(D'_{q-1}:\Omega_{Z,j-2}^{n-2}\to\Omega_{Z,j}^{n-1})=
\Ker(D'_q:\Omega_{Z,j}^{n-1}\to\Omega_{Z,j+2}^n)$ if $q\ne 1$.
}

\medskip\noindent
{\it Proof.}
Let $\phi\in\Omega_{Z,j+2}^n$. There is $\eta\in\Omega_{Z,j}^{n-1}$
such that $\d h\wedge\eta=\phi$ since $j+2>n$.
Then (4.1.2) and (4.1.4) imply
$$
[D'_q\eta]=t\rd_t[\phi]-q[\phi]=(j/2-q)[\phi].
$$
So, replacing $\phi$ with $\phi-\alpha D'_q\eta$
where $\alpha=(j/2-q)^{-1}$, we may assume $[\phi]=0$, i.e.
$\phi\in\d h\wedge\d\Omega_Z^{n-2}$.
Take $\sigma\in \Omega_{Z,j}^{n-2}$ such that $\phi=\d h\wedge\d\sigma$.
Then $D'_q(-q^{-1}\d\sigma)=\phi$, and the assertion~(i) follows.

For the assertion~(ii), let $\eta\in\Omega_{Z,j}^{n-1}$ such that
$D'_q\eta=0$.
Set
$$
\phi=q\d h\wedge\eta=h\d\eta.
$$
Then $t\rd_t[\phi]=q[\phi]$ by (4.1.2), and hence $[\phi]=0$ by
(4.1.4) (using $j\ne 2q$).
Since $H^nA_h^{\ssbull}$ is $t$-torsion-free, we have $[\d\eta]=0$,
i.e. $\d\eta=\d h\wedge\d\sigma$ with $\sigma\in\Omega_{Z,j}^{n-2}$.
Then
$$
\d(D'_{q-1}\sigma)=q\d h\wedge\d\sigma=q\d\eta,
$$
and replacing $\eta$ by $\eta-q^{-1}D'_{q-1}(\sigma)$, we may assume
$d\eta=0$ and hence $\d h\wedge\eta=0$.

In the case $n>2$, this together with (4.1.3) and (4.1.1) implies
$$
\eta=\d\sigma'=-\d h\wedge\d\sigma''\quad\text{with}\,\,\,
\sigma'=\d h\wedge\sigma''\in A_{h,j}^{n-2},\,\,
\sigma''\in\Omega_{h,j-2}^{n-3},
$$
and hence $\eta=(q-1)^{-1}D'_{q-1}(\d\sigma'')$.
So the assertion follows.

In the case $n=2$, we have by (4.1.5) $\eta=\beta h^i\d h$ with
$\beta\in\bC$ if $j$ is even and positive (where $j=2i+2$),
and $\eta=0$ otherwise. If $j=2i+2$, we have
$$
D'_{q-1}h^i=(i-q+1)h^i\d h,
$$
and $i-q+1\ne 0$ by $j\ne 2q$. So the assertion follows.

\medskip\noindent
{\bf 4.3.~Proof of Theorem~2.}
Let $q=n-p$, $i=q-m+1$, $k=(q+1)d$, where $q>m$.
By (2.6) and (2.7)(i), it is enough to treat the case
condition~$(A)$ is satisfied.
With the notation of (2.3), consider the commutative diagram
$$
\CD
0 @>>> (I^{(i)}\Omega^{n})_{k-d} @>>> (I^{(i-1)}\Omega^{n})_{k-d} @>>>
\frac{\displaystyle (I^{(i-1)}\Omega^{n})_{k-d}}
{\displaystyle (I^{(i)}\Omega^{n})_{k-d}} @>>> 0\\
@. @VV{\psi'_a}V @VV{\psi_a}V @VV{\psi''_a}V\\
0 @>>>
\frac{\displaystyle (I^{(i+1)}\Omega^{n+1})_k}
{\displaystyle (fI^{(i-1)}\Omega^{n+1})_k}@>>>
\frac{\displaystyle (I^{(i)}\Omega^{n+1})_k}
{\displaystyle (fI^{(i-1)}\Omega^{n+1})_k}@>>>
\frac{\displaystyle (I^{(i)}\Omega^{n+1})_k}
{\displaystyle (I^{(i+1)}\Omega^{n+1})_k}@>>> 0
\endCD
$$
where $\psi'_a,\psi_a,\psi''_a$ are induced by
$$
\hbox{$D_q:=f\d-q\d f\wedge$
if $a=1$, and $\d f\wedge$ if $a=2$.}
$$
Note that $D_q$ is closely related to (2.1.3).
Using coordinates $x_0,\dots,x_n$, we have
$$
(I^{(i-1)}\Omega^{n})_{k-d}=\mopl_{j=0}^n\,I^{(i-1)}_{k-n-d}\,
\omega_j,\,\,\, (I^{(i)}\Omega^{n+1})_k=I^{(i)}_{k-n-1}\,\omega,\,
\text{etc.,}
$$
where
$\omega=\d x_0\wedge\cdots\wedge\d x_n$ and
$\omega_j=\d x_0\wedge\cdots\widehat{\d x_j}\cdots\wedge\d x_n$.
Then we get
$$
\Coker\,\psi_1=\Gr_F^pH^n(U,\bC),\quad
\Coker\,\psi_2=(I^i/JI^{i-1})_k.
\leqno(4.3.1)
$$
Indeed, the first isomorphism of (4.3.1) follows from (2.1) together
with (2.8), and the second is trivial because $f\in J$.
Note that the assumption of (2.8) is satisfied by condition~$(A)$
for $(k,i)=(qd-n-1,q-m)$, because $i'=q-m$, $k'=qd$ in (2.3.4) with
$j=n-1$.

Since $\rd_jI^{(i)}\subset I^{(i-1)}$, we see that $f\d$ in $\psi'_1$
vanishes, and hence
$$
\Coker\,\psi'_1=\Coker\,\psi'_2.
$$
We will show that $\psi''_a$ is surjective for $a=1,2$
by identifying it with
$$
\mopl_{y\in\Sing Y}\mopl_{j=0}^n\,(\fm_{X,y}^{i-1}/\fm_{X,y}^{i})\,
\omega^{(y)}_j\to\mopl_{y\in\Sing Y}\,(\fm_{X,y}^{i}/\fm_{X,y}^{i+1})\,
\omega^{(y)},
\leqno(4.3.2)
$$
where $\omega^{(y)}_j$, $\omega^{(y)}$ are associated to some
coordinates $x^{(y)}_0,\dots,x^{(y)}_n$ depending on $y\in\Sing Y$.
By the snake lemma, we get then an exact sequence
$$
\Ker\,\psi''_a\buildrel{\rho_a}\over\to\Coker\,\psi'_a\to
\Coker\,\psi_a\to 0.
\leqno(4.3.2)
$$
For $a=1$ this implies the last isomorphism of the formula in Theorem~2.
For the first isomorphism of the formula, we will further show
$$
\Im\,\rho_1=\Im\,\rho_2.
\leqno(4.3.4)
$$

We start with the proof of the surjectivity of $\psi''_a$.
For each $y\in\Sing Y$, we choose appropriate coordinates
$x^{(y)}_0,\dots,x^{(y)}_n$ such that $y$ is given by $(1,0,\dots,0)$,
and
$$
h(z^{(y)}_1,\dots,z^{(y)}_n):=f/(x^{(y)}_0)^d=\msum_{j=1}^n\,
(z^{(y)}_j)^2+\text{higher terms},
$$
where $z^{(y)}_j=x^{(y)}_j/x^{(y)}_0$.
(The last condition is satisfied by using a linear transformation of
$z^{(y)}_1,\dots,z^{(y)}_n$.)
We trivialize $\cO_{X,y}(1)$ by using $x^{(y)}_0$.
Then $\gamma_k^{(i)}$ in (2.3.2) is induced by substituting
$x^{(y)}_0=1$ and $x^{(y)}_j=z^{(y)}_j$ for $j>0$.
So  $z^{(y)}_j$ is identified with $x^{(y)}_j/x^{(y)}_0$.
Since
$$
f_0/(x^{(y)}_0)^{d-1}\in\fm_{X,y}^2,\quad
f_j/(x^{(y)}_0)^{d-1}=2z^{(y)}_j\,\,\,\text{in}\,\,\,
\fm_{X,y}/\fm_{X,y}^2\,(j\ne 0),
$$
we see that $\psi''_a$ is identified with (4.3.2).
Indeed, the assertion is equivalent to the surjectivity of
$\gamma_{k-n-d}^{(i-1)}$ and $\gamma_{k-n-1}^{(i)}$ in (2.3.2).
But the first surjectivity follows from condition~$(A)$ for $(k,i)=
(qd-n,q-m+1)$, and the second is reduced to the first by using a
commutative diagram as above together with the surjectivity of the
morphism~(4.3.2) induced by $df\wedge$, see (4.1.1).
So $\psi''_a$ is identified with (4.3.2), and we get also the
surjectivity of $\psi''_2$.

The morphisms (4.3.2) induced by $D_q$ and $df\wedge$ are compatible
with the direct sum over $y\in\Sing Y$
(using the pull-back by the surjection $\gamma_{k-n-d}^{(i-1),(y)}$,
see (2.3.3)), and moreover, the restriction of $\psi''_1$
$$
\mopl_{j=1}^n\,(\fm_{X,y}^{i-1}/\fm_{X,y}^{i})\,\omega^{(y)}_j\to
(\fm_{X,y}^{i}/\fm_{X,y}^{i+1})\,\omega^{(y)}
\leqno(4.3.5)
$$
is identified with $D'_q$ in Lemma~(4.2).
Here $j:=i+n-2\ne 2q$ since $q>m$.
So $\psi''_1$ is also surjective, and the kernel of (4.3.5)
does not contribute to $\Im\,\rho_1$ by using
$$
D_{q-1}:(I^{i-2}\Omega^{n}/I^{i-1}\Omega^{n})_{k-2d}\to
(I^{i-1}\Omega^{n}/I^i\Omega^{n})_{k-d},
$$
because it lifts $D'_{q-1}$ in Lemma~(4.2) and satisfies
$D_q\scirc D_{q-1}=0$.
For $a=2$, the kernel of (4.3.5) induced by $\d f\wedge$
does not contribute to $\Im\,\rho_2$ by a similar argument
using (4.1.1).

Thus it is enough to consider the contribution to $\Im\,\rho_a$ of
$(\fm_{X,y}^{i-1}/\fm_{X,y}^{i})\,\omega_0^{(y)}$ which is contained
in the kernel of (4.3.2) for $a=1,2$.
Since $\rd/\rd x^{(y)}_0$ preserves the maximal ideal of $R$ generated by
$x_j^{(y)}\,(j\ne 0)$,
it does not contribute to $\Im\,\rho_1$ using the pull-back by the
surjection
$\gamma_{k-n-1}^{(i),(y)}$ in (2.3.3).
Then the contributions to $\Im\,\rho_a$ for $a=1,2$ are both given by
using the pull-back by the surjection $\gamma_{k-n-1}^{(i),(y)}$
together with the multiplication by $\rd f/\rd x^{(y)}_0$.
So we get (4.3.4).
(Note that the assertion (4.3.4) is independent of the choice of
the coordinates, and the obtained isomorphism in the formula of
Theorem~2 is well defined.)
This completes the proof of Theorem~2.

\medskip\noindent
{\bf 4.4.~Proof of Corollary~1.}
Let $q=n-p$, $i=q-m+1$, $k=(q+1)d$.
Since $q>m\ge 1$, the condition in (2.5.1) in the case $d\ge 3$ is
satisfied when $k,i$ in (2.5.1) are $k-n-1,i$ or $k-n-d,i-1$
(i.e. we have $k-n-1\ge 2i$ and $k-n-d\ge 2i-2$),
see (2.4.3) for the case $d=2$.
Moreover, $m(d-1)-p>0$ in the case $d\ge 3$,
because $n-p>m\ge 1$.
If $d=2$, then $|\Sing\,Y|=1$ and (2.3.5) is satisfied for $e=0$.
So the assertion follows from Theorem~2 and Lemma~(2.5).

\medskip
From Theorem~2 we can deduce the following.

\medskip\noindent
{\bf 4.5.~Corollary.} {\it Let $Y_s$ be an equisingular family of
hypersurfaces in $\bP^n$ which are parametrized by a smooth variety $S$ and
whose singularities are ordinary double points.
Assume condition $(A)$ for $q$ in Theorem $2$ is satisfied for any
$s\in S$ if $q>m$, and assume the same with $q$ replaced by $q-1$
if $q-1>m$.
Set $U_s=\bP^n\setminus Y_s$.
Then for a vector field $\theta$ on S, we have a commutative diagram
$$
\CD
\Gr_F^{n-q+1}H^n(U_s,\bP) @>{\Gr_F\nabla_{\theta}}>>
\Gr_F^{n-q}H^n(U_s,\bP)\\
@| @|\\
(I_s^{(q-m)}/I_s^{(q-m-1)}J_s)_{qd-n-1} @>{-q(\theta f)_s}>>
(I_s^{(q-m+1)}/I_s^{(q-m)}J_s)_{(q+1)d-n-1},\\
\endCD
$$
where the vertical isomorphisms are given by Theorems~{\rm 2} and
{\rm (2.2)}.
}
 
\medskip\noindent
{\it Proof.}
The action of $\theta$ on the relative de Rham cohomology
can be calculated by $\iota_{\theta}\scirc\d$.
So the assertion follows from Theorem~2 (using
$\iota_{\xi}\scirc\d+\d\scirc\iota_{\xi}=L_{\xi}$ and
$\iota_{\theta}\scirc\iota_{\xi}=-\iota_{\xi}\scirc\iota_{\theta}$),
since the cohomology class is represented by the middle term of (3)
in Theorem~2.

\medskip\noindent
{\bf 4.6.~Remarks.}
(i) By Varchenko~[22] (conjectured by Arnold) and [13], we have
$$
\aligned
|\Sing Y|&\le \msum_{(n-2)/2+1<i\le nd/2}\,C(n,d,i)=C(n+1,d,[nd/2]+1)\\
&=\msum_{i\ge 0}\hbox{$\binom{n+1}{i}\binom{[nd/2]-i(d-1)}{n}$}<
\hbox{$\binom{[nd/2]}{n}$},
\endaligned
\leqno(4.6.1)
$$
where $C(n,d,i)$ is as in (3.3.1).
(This also follows from (3.5.2) applied to a
hypersurface in $\bP^{n+1}$ or $\bP^{n+2}$ defined by
$f+x_{n+1}^d$ or $f+x_{n+1}^d+x_{n+2}^d$.)

If $n$ is even (i.e. if $n=2m$), then (4.6.1) implies
$$
\hbox{$\binom{e+n}{n}\ge\binom{md+1}{n}>\binom{[nd/2]}{n}>
|\Sing Y|\,\,$ for
$q\ge m+1$,}
\leqno(4.6.2)
$$
and it is possible that condition~$(B)$ in Theorem~2 is satisfied.
However, if $n$ is odd (i.e. $n=2m+1$), then (4.6.2) does not hold
and condition~$(B)$ cannot be satisfied,
for example, if $n=3$, $q=2$, and $Y$ is a Kummer surface with 16
ordinary double points where $d=4$, $e=2$ and $\binom{e+n}{n}=10$
(see (4.7)(ii) below), or $Y$ has $65$ ordinary double points with
$d=6$ in [1] where $e=4$ and $\binom{e+n}{n}=35$.

(ii) As for condition~$(A)$, we get matrices of size $(M,N)$ with
$$
\hbox{$M=\binom{i-1+n}{n}|\Sing Y|,\quad N=\binom{k+n}{n}$},
$$
where $(k,i)=(qd-n,q-m+1)$ and $(qd-n-1,q-m)$, see Remark (2.7)(ii).

\medskip\noindent
{\bf 4.7.~Examples.}
(i) Assume $n=3$, $d=4$, and
$$
f=\msum_{i=0}^3\,x_i^4-\msum_{0\le i<j\le 3}\,2x_i^2x_j^2
\quad\text{so that}\quad
f_j=4x_j(x_j^2-\msum_{i\ne j}x_i^2).
$$
This has 12 ordinary double points.
Indeed, there are 2 singular points defined by
$x_i^2=x_j^2$ and $x_k^2=0 \,(k\in[0,3]\setminus\{i,j\})$ for each
$\{i,j\}\subset[0,3]$ with $i\ne j$.
So condition~$(B)$ cannot be satisfied for $q=2$
since $12>\binom{e+n}{n}=10$ in the notation of Remark~(4.6).
However, condition~$(A)$ seems to be satisfied for $q=2$ where
$(M,N)=(48,56)$ and $(12,35)$ in the notation of Remark~(4.6)(ii).

(ii) Assume $Y$ is a singular Kummer surface defined by
$$
f=\msum_{i=0}^3\,x_i^4-\msum_{0\le i<j\le 3}\,x_i^2x_j^2
\quad\text{so that}\quad
f_j=2x_j(2x_j^2-\msum_{i\ne j}x_i^2).
$$
This has 16 ordinary double points.
Indeed, there are 4 singular points defined by $x_k=0$ and
$x_i^2=1 \,(i\ne k)$ for each $k=0,\dots,3$.
Condition~$(A)$ for $q=2$ cannot be satisfied since
$(M,N)$ can be $(64,56)$ in the notation of Remark~(4.6)(ii).
However, it seems that
$$
\dim\,I^{(2)}_8= \dim\,(I^2)_8=\hbox{$\binom{8+3}{3}$}-4\,
|\Sing Y|=101,\,\,\,
\dim\,(IJ)_8=100,
$$
so that at least a noncanonical isomorphism still holds in Conjecture~1
for $q=2$.

\medskip
(iii) Since it is not easy to calculate the right-hand side of
Conjecture~1 for the Barth surface~[1], we consider the case $Y$ is
defined by
$$
f=(\msum_{i=0}^3\,x_i^2)^3-\msum_{i=0}^3\,x_i^6
\quad\text{so that}\quad
f_j=6x_j((\msum_{i=0}^3\,x_i^2)^2-x_j^4).
$$
This has 52 ordinary double points.
Indeed, there are 4 singular points defined by $x_i=1$ and
$x_k=0 \,(k\ne i)$ for $i=0,\dots,3$,
and there are 4 singular points defined by $x_i=1$, $x_j=0$ and
$x_k^2=-1 \,(k\in[0,3]\setminus\{i,j\})$ for each $(i,j)\in[0,3]^2
\setminus\{\hbox{diagonal}\}$.
Condition~$(B)$ cannot be satisfied since $\binom{e+n}{n}=35<52$,
but it is not clear whether condition~$(A)$ is satisfied where
$(M,N)=(208,220)$ and $(52,165)$.
It seems that
$$
\dim\,I^{(2)}_{14}= \dim\,(I^2)_{14}=\hbox{$\binom{14+3}{3}$}-4\,
|\Sing Y|=472,\,\,\,
\dim\,(IJ)_{14}=462,
$$
so that at least a noncanonical isomorphism still holds in Conjecture~1
for $q=2$.

\bigskip
\ver

\end{document}